\documentclass[11pt]{article}
\usepackage{epic,latexsym,amssymb}
\usepackage{color}
\usepackage{tikz}
\usepackage{amsfonts,epsf,amsmath}
\usepackage{cite}

\newtheorem{thm}{Theorem}
\newtheorem{lem}{Lemma}

\newtheorem{claim}{Claim}
\newtheorem{subclaim}{Claim}[claim]

\newcommand{\qed}{$\Box$}

\newcommand{\smallqed}{{\tiny ($\Box$)}}

\newcommand{\ii}{\iota}

\newcommand{\proof}{\noindent\textbf{Proof. }}

\let\oldenumerate\enumerate
\renewcommand{\enumerate}{
  \oldenumerate
  \setlength{\itemsep}{0pt}
  \setlength{\parskip}{0pt}
  \setlength{\parsep}{0pt}
}

\begin{document}

\title{Partial domination of maximal outerplanar graphs}

\author{Peter Borg$^1$ \, and \, Pawaton Kaemawichanurat$^2$\thanks{Research supported in part by Skill Development Grant 2018, King Mongkut's University of Technology Thonburi.}
\\ \\
$^1$Department of Mathematics \\
Faculty of Science \\
University of Malta\\
Malta\\
$^2$Theoretical and Computational Science Center \\
and Department of Mathematics, Faculty of Science,\\
King Mongkut's University of Technology Thonburi \\
Bangkok, Thailand\\
\small \tt Email: $^{1}$peter.borg@um.edu.mt, $^{2}$pawaton.kae@kmutt.ac.th}

\date{}
\maketitle

\begin{abstract}
Several domination results have been obtained for maximal outerplanar graphs (mops). The classical domination problem is to minimize the size of a set $S$ of vertices of an $n$-vertex graph $G$ such that $G - N[S]$, the graph obtained by deleting the closed neighborhood of $S$, is null. A classical result of Chv\'{a}tal is that the minimum size is at most $n/3$ if $G$ is a mop. Here we consider a modification by allowing $G - N[S]$ to have isolated vertices and isolated edges only. Let $\iota_1(G)$ denote the size of a smallest set $S$ for which this is achieved. We show that if $G$ is a mop on $n \geq 5$ vertices, then $\iota_{1}(G) \leq n/5$. We also show that if $n_2$ is the number of vertices of degree $2$, then $\iota_{1}(G) \leq \frac{n+n_2}{6}$ if $n_2 \leq \frac{n}{3}$, and $\iota_1(G) \leq \frac{n-n_2}{3}$ otherwise. We show that these bounds are best possible.
\end{abstract}


\section{Introduction}
Let $G$ be a simple graph with vertex set $V(G)$ and edge set $E(G)$. We denote the \emph{degree} of $v$ in $G$ by $d_G(v)$. The \emph{open neighborhood} $N_G(v)$ of a vertex $v$ of $G$ is the set of neighbors of $v$, that is, $N_G(v) = \{u \in V(G) \colon uv \in E(G)\}$. The \emph{closed neighborhood} of $v$ is $N_G[v] = N_G(v) \cup \{v\}$. The \emph{neighborhood} of a subset $S$ of $V(G)$ is the set $N_G(S) = \cup_{v \in S} N_G(v)$. The \emph{closed neighborhood} of a subset $S$ of $V(G)$ is the set $N_G[S] = N_G(S) \cup S$. For $S \subseteq V(G)$, the subgraph of $G$ induced by $S$ is denoted by $G[S]$. The subgraph obtained from $G$ by deleting all vertices in $S$ (and all edges incident to vertices in $S$) is denoted by $G - S$.

A subset $S$ of $V(G)$ is a \emph{dominating set} of $G$ if each vertex in $V(G) \setminus S$ is adjacent to at least one vertex in $S$. The \emph{domination number of $G$} is the size of a smallest dominating set of $G$ and is denoted by $\gamma(G)$. Given a graph $H$, a subset $S$ of $V(G)$ is an \emph{$H$-isolating set of $G$} if $G - N_{G}[S]$ does not contain a copy of $H$. Obviously, $S$ is a dominating set if $H = K_{1}$. The \emph{$H$-isolation number of $G$} is the size of a smallest $H$-isolating set of $G$ and is denoted by $\ii(G, H)$. If $H = K_{1, k + 1}$ for some $k \geq 0$, then we may abbreviate $\ii(G, H)$ to $\ii_{k}(G)$. The study of isolating sets is an appealing and natural direction in domination theory that was introduced recently by Caro and Hansberg~\cite{CaHa17}.

A \emph{triangulated disc} is a simple planar graph whose interior faces are triangles. A \emph{maximal outerplanar graph} is a triangulated disc whose exterior face (the unbounded face) contains all vertices. Hence, a maximal outerplanar graph can be embedded in the plane such that all vertices lie on the boundary of the exterior face and all interior faces are triangles. We abbreviate the term \emph{maximal outerplanar graph} to \emph{mop}. O'Rourke~\cite{Ro87} pointed out that every mop has a unique Hamiltonian cycle. Thus, the Hamiltonian cycle of a mop is the boundary of the mop. This paper's notation and terminology on mops follows that of~\cite{LeZuZy17}; in particular, an edge belonging to the Hamiltonian cycle of a mop is called a \emph{Hamiltonian edge}, while any other edge of the mop is called a \emph{diagonal}. A \emph{fan} of order $n \ge 3$, denoted $F_n$, is the mop obtained from a path $P_{n-1}$ by adding a new vertex $v$ and joining it to every vertex of the path. We say that $v$ is the \emph{center} of $F_n$.

Domination in mops has been extensively studied since 1975. In the classical paper \cite{Ch75}, Chv\'{a}tal essentially proved that the domination number of an $n$-vertex mop is at most $n/3$. Other proofs were obtained by Fish~\cite{Fi78} and by Matheson and Tarjan~\cite{MaTa96}. Caro and Hansberg~\cite{CaHa17} proved that the $K_{1, 1}$-isolation number of a mop of order $n \ge 4$ is at most $n/4$.

\begin{thm}
\label{t:known1}
If $G$ is a mop of order $n$, then the following hold: \\
\indent {\rm (a)} {\rm \cite{Ch75,Fi78,MaTa96}} If $n \ge 3$, then $\gamma(G) \le \frac{n}{3}$. \\
\indent {\rm (b)} {\rm \cite{CaHa17}} If $n \ge 4$, then $\ii_0(G) \le \frac{n}{4}$.
\end{thm}

For results on other types of domination in mops, we refer the reader to \cite{CaHa17,CaWa13,CaCaHeMa16,DoHaJo16,DoHaJo17,HeKa18,KaJi,To13}.

Consider any $n$-vertex mop $G$. Theorem \ref{t:known1}(a) tells us that there exists some $S \subseteq V(G)$ such that $|S| \leq n/3$ and $G-N_G[S]$ is null. Theorem \ref{t:known1}(b) tells us that there exists some $S \subseteq V(G)$  such that $|S| \leq n/4$ and $G-N_G[S]$ has isolated vertices only. In this paper, we establish sharp upper bounds on the size of a smallest subset $S$ of $V(G)$ such that $G-N_G[S]$ has isolated vertices and isolated edges only, that is, the edges of $G-N_G[S]$ form a matching (that is, they are pairwise disjoint). The bounds are detailed in the following two theorems and the proofs are given in Section~\ref{S:main1}.

\begin{thm}
\label{t:main1}
If $G$ is a mop of order $n \ge 5$, then \[\ii_{1}(G) \leq \frac{n}{5}.\]
\end{thm}

\begin{thm}
\label{t:main2}
If $G$ is a mop of order $n \ge 5$ with $n_2$ vertices of degree~$2$, then
\[ \ii_1(G)  \le \left\{\begin{array}{ll}
                     \frac{n + n_2}{6} & \mbox{if } n_2 \le \frac{n}{3},\\
                     \\
                     \frac{n - n_2}{3} & otherwise.
               \end{array}\right. \]
\end{thm}

The bounds in Theorems~\ref{t:main1} and \ref{t:main2} are sharp. For example, let $F^{1}_{5},F^{2}_{5}, \dots, F^{t}_{5}$ be $t \ge 2$ vertex disjoint fans of order~$5$. From each $F^{i}_{5}$, we choose a vertex of degree~$2$ and its neighbor of degree~$3$. Then, we join these $2t$ vertices by edges to form a mop. The resulting graph $G_t$ is a mop of order $n = 5t$ with $n_2 = t$ vertices of degree~$2$. An example of the graph $G_t$ with $t = 4$ is illustrated in Figure 1.

\setlength{\unitlength}{0.6cm}
\begin{center}
\begin{picture}(15, 6)

\put(2, 5){\circle*{0.2}}
\put(2, 5){\line(1, -1){1}}
\put(2, 5){\line(-1, -1){1}}
\put(1, 2){\circle*{0.2}}
\put(3, 2){\circle*{0.2}}
\put(1, 4){\circle*{0.2}}
\put(3, 4){\circle*{0.2}}
\put(1, 4){\line(1, 0){2}}
\put(1, 4){\line(1, 0){2}}
\put(1, 2){\line(0, 1){2}}
\put(3, 2){\line(0, 1){2}}
\put(3, 2){\line(-1, 1){2}}
\put(1, 2){\line(1, 0){14}}

\put(6, 5){\circle*{0.2}}
\put(6, 5){\line(1, -1){1}}
\put(6, 5){\line(-1, -1){1}}
\put(5, 2){\circle*{0.2}}
\put(7, 2){\circle*{0.2}}
\put(5, 4){\circle*{0.2}}
\put(7, 4){\circle*{0.2}}
\put(5, 4){\line(1, 0){2}}
\put(5, 4){\line(1, 0){2}}
\put(5, 2){\line(0, 1){2}}
\put(7, 2){\line(0, 1){2}}
\put(7, 2){\line(-1, 1){2}}

\put(10, 5){\circle*{0.2}}
\put(10, 5){\line(1, -1){1}}
\put(10, 5){\line(-1, -1){1}}
\put(9, 2){\circle*{0.2}}
\put(11, 2){\circle*{0.2}}
\put(9, 4){\circle*{0.2}}
\put(11, 4){\circle*{0.2}}
\put(9, 4){\line(1, 0){2}}
\put(9, 4){\line(1, 0){2}}
\put(9, 2){\line(0, 1){2}}
\put(11, 2){\line(0, 1){2}}
\put(11, 2){\line(-1, 1){2}}

\put(14, 5){\circle*{0.2}}
\put(14, 5){\line(1, -1){1}}
\put(14, 5){\line(-1, -1){1}}
\put(13, 2){\circle*{0.2}}
\put(15, 2){\circle*{0.2}}
\put(13, 4){\circle*{0.2}}
\put(15, 4){\circle*{0.2}}
\put(13, 4){\line(1, 0){2}}
\put(13, 4){\line(1, 0){2}}
\put(13, 2){\line(0, 1){2}}
\put(15, 2){\line(0, 1){2}}
\put(15, 2){\line(-1, 1){2}}

\qbezier(1, 2)(6, -3.5)(11, 2)
\qbezier(1, 2)(5, -2)(9, 2)
\qbezier(1, 2)(4, -1)(7, 2)
\qbezier(3, 2)(5, -0.1)(7, 2)
\qbezier(1, 2)(7, -5)(15, 2)
\qbezier(11, 2)(13, -0)(15, 2)

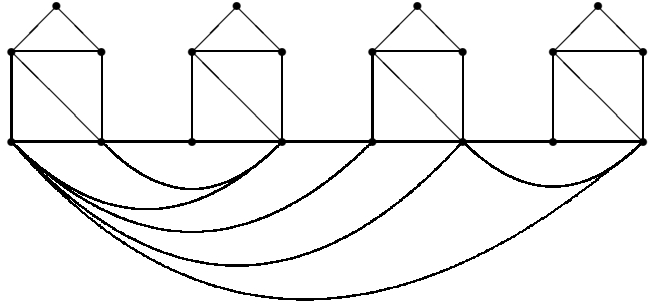
\put(0, -3){\footnotesize{\textbf{\textbf{Figure 1 :}} The mop $G_{4}$, which satisfies $\ii_1(G_{4}) = n/5 = 20/5 = 4$}}
\put(3, -4){\footnotesize{and $\ii_1(G_{4}) = (n+n_2)/6 = (20 + 4)/6 = 4$.}}
\end{picture}
\end{center}
\vskip 70 pt

\noindent Clearly, for a $K_{1, 2}$-isolating set $S$ of $G_t$, we have $|S \cap V(F^{i}_{5})| \geq 1$ because $G_t - N_{G_{t}}[S]$ does not contain $K_{1, 2}$. Thus, $\ii_1(G_t) \ge t$. However, $G_k$ has a $K_{1, 2}$-isolating set containing the vertex of degree~$4$ of each $F^{i}_{5}$. Consequently, $\ii_1(G_t) = t = (5t)/5 = n/5$ and $\ii_1(G_t) = t = (5t + t)/6 = (n + n_2)/6$.

\setlength{\unitlength}{0.6cm}
\begin{center}
\begin{picture}(15, 5)
\put(1, 2){\circle*{0.2}}
\put(2, 4){\circle*{0.2}}
\put(3, 2){\circle*{0.2}}
\put(3, 2){\circle*{0.2}}
\put(1, 4){\circle*{0.2}}
\put(3, 4){\circle*{0.2}}
\put(1, 4){\line(1, 0){2}}
\put(1, 4){\line(1, 0){2}}
\put(1, 2){\line(0, 1){2}}
\put(3, 2){\line(0, 1){2}}
\put(3, 2){\line(-1, 2){1}}
\put(1, 2){\line(1, 2){1}}
\put(1, 2){\line(1, 0){14}}

\put(5, 2){\circle*{0.2}}
\put(6, 4){\circle*{0.2}}
\put(7, 2){\circle*{0.2}}
\put(5, 4){\circle*{0.2}}
\put(7, 4){\circle*{0.2}}
\put(5, 4){\line(1, 0){2}}
\put(5, 4){\line(1, 0){2}}
\put(5, 2){\line(0, 1){2}}
\put(7, 2){\line(0, 1){2}}
\put(7, 2){\line(-1, 2){1}}
\put(5, 2){\line(1, 2){1}}

\put(9, 2){\circle*{0.2}}
\put(10, 4){\circle*{0.2}}
\put(11, 2){\circle*{0.2}}
\put(9, 4){\circle*{0.2}}
\put(11, 4){\circle*{0.2}}
\put(9, 4){\line(1, 0){2}}
\put(9, 4){\line(1, 0){2}}
\put(9, 2){\line(0, 1){2}}
\put(11, 2){\line(0, 1){2}}
\put(11, 2){\line(-1, 2){1}}
\put(9, 2){\line(1, 2){1}}

\put(13, 2){\circle*{0.2}}
\put(14, 4){\circle*{0.2}}
\put(15, 2){\circle*{0.2}}
\put(13, 4){\circle*{0.2}}
\put(15, 4){\circle*{0.2}}
\put(13, 4){\line(1, 0){2}}
\put(13, 4){\line(1, 0){2}}
\put(13, 2){\line(0, 1){2}}
\put(15, 2){\line(0, 1){2}}
\put(15, 2){\line(-1, 2){1}}
\put(13, 2){\line(1, 2){1}}

\qbezier(1, 2)(6, -3.5)(15, 2)
\qbezier(1, 2)(5, -2)(9, 2)
\qbezier(3, 2)(6, -1)(9, 2)
\qbezier(3, 2)(5, -0.1)(7, 2)
\qbezier(9, 2)(12, 0)(15, 2)
\qbezier(9, 2)(11, 1)(13, 2)

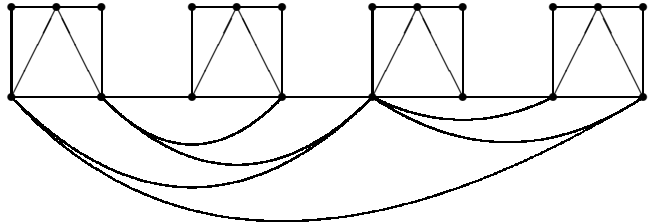
\put(-1, -2){\footnotesize{\textbf{\textbf{Figure 2 :}} The mop $H_{4}$, which satisfies $\ii_1(H_{4}) = (n - n_2)/3 = (20 - 8)/3 = 4$.}}
\end{picture}
\end{center}
\vskip 40 pt

\indent To see the sharpness of the bound $\ii_1(G) \le \frac{n - n_2}{3}$, we first let  $F^{1}_{5}, F^{2}_{5}, \dots, F^{t}_{5}$ be $t \ge 2$ fans of order~$5$. We choose the two vertices of degree~$3$ from each $F^{i}_{5}$. Then we join these $2t$ vertices by edges to form a mop. The resulting graph $H_t$ is a mop of order~$n = 5t$ with $n_2 = 2t$ vertices of degree~$2$. An example of the graph $H_t$ with $t = 4$ is illustrated in Figure 2. Clearly, for a $K_{1, 2}$-isolating set $S$ of $H_t$, we have $|S \cap V(F^{i}_{5})| \geq 1$ because $H_t - N_{H_{t}}[S]$ does not contain $K_{1, 2}$. Thus, $\ii_1(H_t) \ge t$. However, $H_t$ has a $K_{1, 2}$-isolating set containing the vertex of degree~$4$ of each $F^{i}_{5}$. Consequently, $\ii_1(H_t) = t = (5t - 2t)/3 = (n - n_2)/3$.

We conclude this section by comparing the bound $n/5$ with $(n+n_2)/6$ and $(n-n_2)/3$. We see that $n/5 < (n+n_2)/6$ when $n/5 < n_2$. Moreover, $n/5 < (n-n_2)/3$ when $n_2 < 2n/5$. Thus, the bound in Theorem \ref{t:main1} is better than that in Theorem \ref{t:main2} when $n/5 < n_2 < 2n/5$. We will show that the bound $n/5$ is sharp when $n/5 < n_2 < 2n/5$. Let $F_1$, $F_2$ and $F_3$ be vertex-disjoint fans of order $5$ such that $V(F_1)=\{x_0,x_1,\dots,x_{4}\},V(F_2)=\{y_0,y_1,\dots,y_{4}\}$ and $V(F_3)=\{z_0,z_1,\dots,z_{4}\}$, where $x_0,y_0,z_0$ are the centers of the fans and $x_1,x_{4},y_1,y_{4},z_1,z_{4}$ are the vertices of degree~$2$ of the fans. Let $A_{15}$ be the graph obtained by taking the union of $F_1, F_2, F_3$ and a mop that has $\{z_1,z_2,x_2,x_3,y_2,y_3\}$ as its vertex set and has $z_1z_2$, $x_2x_3$ and $y_2y_3$ among its Hamiltonian edges. Clearly, $A_{15}$ is a mop of order $15$ with $5$ vertices of degree~$2$. This graph is illustrated in Figure 3. Then, we let $A^{1}_{15},A^{2}_{15},\dots,A^{t}_{15}$ be $t$ vertex-disjoint copies of $A_{15}$. The graph $B_t$ is constructed by taking the union of $A^{1}_{15},A^{2}_{15},\dots,A^{t}_{15}$ and a mop that has $\{x^{1}_{3}, y^{1}_{3}, x^{2}_{3}, y^{2}_{3}, \dots, x^{t}_{3}, y^{t}_{3}\}$ as its vertex set and has $x^{1}_{3}y^{1}_{3}, x^{2}_{3}y^{2}_{3}, \dots, x^{t}_{3}y^{t}_{3}$ among its Hamiltonian edges.
Clearly, $B_t$ is a mop of order $n=15t$ with $n_2=5t$ vertices of degree~$2$. Thus, $n/5 <n_2=5n/15 < 2n/5$. Clearly, each $K_{1,2}$-isolating set of $B_t$ needs at least one vertex from each fan $F^{i}_{j}$.
Thus, $\ii_k(B_t) \ge 3t$. Moreover, $\{x^{i}_{0} \colon 1 \le i \le t\} \cup \{y^{i}_{0} \colon 1 \le i \le t\} \cup \{z^{i}_{0} \colon 1 \le i \le t\}$ is a $K_{1,2}$-isolating set of $B$. Thus, $\iota_1(B_t) = 3t = n/5$.

\setlength{\unitlength}{0.9cm}
\begin{center}
\begin{picture}(1, 8)

\put(0, 0){\circle*{0.2}}
\put(0, -0.25){\footnotesize{$x_{3}$}}

\put(0, 2){\circle*{0.2}}
\put(0.25, 1.75){\footnotesize{$x_{2}$}}

\put(-2, 2){\circle*{0.2}}
\put(-2.4, 2.25){\footnotesize{$x_{0}$}}

\put(-4, 2){\circle*{0.2}}
\put(-4.25, 2.25){\footnotesize{$x_{4}$}}

\put(-2, 4){\circle*{0.2}}
\put(-2, 4.25){\footnotesize{$x_{1}$}}

\put(2, 0){\circle*{0.2}}
\put(2.25, -0.25){\footnotesize{$y_{3}$}}

\put(2, 2){\circle*{0.2}}
\put(2.25, 1.75){\footnotesize{$y_{2}$}}

\put(4, 4){\circle*{0.2}}
\put(4, 4.25){\footnotesize{$y_{1}$}}

\put(4, 2){\circle*{0.2}}
\put(4.25, 2.25){\footnotesize{$y_{0}$}}

\put(6, 2){\circle*{0.2}}
\put(6.25, 1.75){\footnotesize{$y_{4}$}}

\put(0, 4){\circle*{0.2}}
\put(0.25, 3.75){\footnotesize{$z_{1}$}}

\put(2, 4){\circle*{0.2}}
\put(2.25, 3.75){\footnotesize{$z_{2}$}}

\put(2, 6){\circle*{0.2}}
\put(2.25, 6){\footnotesize{$z_{3}$}}

\put(1, 7){\circle*{0.2}}
\put(1, 7.25){\footnotesize{$z_{4}$}}

\put(0, 6){\circle*{0.2}}
\put(-0.25, 6.25){\footnotesize{$z_{0}$}}

\put(0, 0){\line(0, 2){2}}
\put(0, 2){\line(-2, 0){2}}
\put(-2, 2){\line(-2, 0){2}}
\put(0, 0){\line(-1, 1){2}}
\put(0, 0){\line(-2, 1){4}}
\put(-2, 2){\line(0, 1){2}}
\put(0, 2){\line(-1, 1){2}}
\put(0, 2){\line(0, 1){2}}
\put(0, 4){\line(0, 1){2}}
\put(0, 0){\line(1, 0){2}}
\put(2, 0){\line(0, 1){2}}
\put(0, 2){\line(1, 0){2}}
\put(0, 0){\line(1, 1){2}}
\put(2, 2){\line(0, 1){2}}
\put(2, 4){\line(0, 1){2}}
\put(0, 4){\line(1, 0){2}}
\put(0, 6){\line(1, 0){2}}
\put(0, 6){\line(1, 1){1}}
\put(2, 6){\line(-1, 1){1}}
\put(0, 2){\line(1, 1){2}}
\put(2, 4){\line(-1, 1){2}}
\put(2, 0){\line(1, 1){2}}
\put(2, 2){\line(1, 0){2}}
\put(4, 2){\line(1, 0){2}}
\put(2, 0){\line(2, 1){4}}
\put(2, 2){\line(1, 1){2}}
\put(4, 2){\line(0, 1){2}}

\put(-1, -1){\footnotesize{\textbf{Figure 3 : }The mop \textbf{$A_{15}$}.}}
\end{picture}
\end{center}
\vskip 50 pt

\section{Proofs of the upper bounds}
\label{S:main1}

In this section, we prove Theorems~\ref{t:main1} and \ref{t:main2}. We apply results of O'Rourke~\cite{Ro83} in computational geometry that were used in a new proof by Lema\'{n}ska, Zuazua and Zylinski~\cite{LeZuZy17} of an upper bound by Dorfling, Hattingh and Jonck \cite{DoHaJo16} on the size of a total dominating set~(a set $S$ of vertices such that each vertex of the graph is adjacent to a vertex in $S$) of a mop. Before stating these results, we make a related straightforward observation that we will also use.

Given three mops $G$, $G_1$ and $G_2$, we say that a diagonal $d$ of $G$ \emph{partitions $G$ into $G_1$ and $G_2$} if $G$ is the union of $G_1$ and $G_2$, $V(G_1) \cap V(G_2) = d$ and $E(G_1) \cap E(G_2) = \{d\}$.

\begin{lem}\label{diagonallemma} If $d$ is a diagonal of a mop $G$, then $d$ partitions $G$ into two mops $G_1$ and $G_2$.
\end{lem}

\proof Let $x_1x_2 \dots x_nx_1$ be the Hamiltonian cycle $C$ of $G$. We may assume that $d = x_1x_i$ for some $i \in \{2, 3, \dots, n\}$. Let $C_1$ be the cycle $x_1x_2 \dots x_ix_1$ of $G$, and let $C_2$ be the cycle $x_1x_ix_{i+1} \dots x_nx_1$. For each $i \in \{1, 2\}$, let $G_i$ be the subgraph of $G$ induced by $V(C_i)$. Then, $V(G_1) \cap V(G_2) = d$ and $E(G_1) \cap E(G_2) = \{d\}$. Each face of $G_i$ is a face of $G$ in the interior of $C_i$ and hence a triangle; thus, $G_i$ is a mop. Since $G$ is a mop and $x_1x_i$ is a diagonal of $G$, $G$ has no edge with one vertex in $V(C_1) \backslash \{x_1, x_i\}$ and the other vertex in $V(C_2) \backslash \{x_1, x_i\}$. Thus, $E(G) = E(G_1) \cup E(G_2)$.\hfill{\qed}

The next lemma was proved by Chv\'{a}tal in \cite{Ch75} and is restated in \cite[Lemma 1.1]{Ro87}, and the lemma following it is an extension by O'Rourke \cite{Ro83} that has a central role in our proofs.

\begin{lem}{\rm (\cite{Ch75})}
\label{mop:4-6}
If $G$ is a mop of order $n \geq 8$, then $G$ has a diagonal $d$ that partitions it into two mops $G_1$ and $G_2$ such that $G_1$ has exactly $4$, $5$ or $6$ Hamiltonian edges of $G$.
\end{lem}

\begin{lem}{\rm (\cite{Ro83})}
\label{mop:5-8}
If $G$ is a mop of order $n \geq 10$, then $G$ has a diagonal $d$ that partitions it into two mops $G_1$ and $G_2$ such that $G_1$ has exactly $5$, $6$, $7$ or $8$ Hamiltonian edges of $G$.
\end{lem}

For a graph $G$ and an edge $uv$ of $G$, the \emph{edge contraction} of $G$ along $uv$ is the graph obtained from $G$ by deleting $u$ and $v$ (and all incident edges), adding a new vertex $x$, and making $x$ adjacent to the vertices in $N_G(\{u,v\}) \setminus \{u,v\}$ only. Recall that every mop can be embedded in a plane so that the exterior face contains all vertices. By looking at polygon corners as vertices, we have that a mop is a triangulation of a simple polygon, meaning that its boundary is the polygon and its interior faces are triangles.

\begin{lem}{\rm (\cite{Ro83})}
\label{contraction}
If $G$ is a triangulation of a simple polygon $P$, $e$ is a Hamiltonian edge of $G$, and $G'$ is the edge contraction of $G$ along $e$, then $G'$ is a triangulation of some simple polygon $P'$.
\end{lem}

If $G$ is a graph and $I \subseteq V(G)$ such that $uv \notin E(G)$ for every $u, v \in I$, then $I$ is called an \emph{independent set of $G$}.

For the proof of Theorem \ref{t:main2}, we establish the following result about vertices of degree $2$.

\begin{lem}\label{maxofn2}
If $G$ is a mop of order $n \geq 4$, then the set of vertices of $G$ of degree $2$ is an independent set of $G$ of size at most $\frac{n}{2}$.
\end{lem}

\proof Let $V_2$ be the set of vertices of $G$ of degree~$2$. Let $x_0x_1 \ldots x_{n-1} x_0$ be the unique Hamiltonian cycle of $G$ and hence the boundary of the exterior face of $G$. For each $i \in \{0, 1, \dots, n-1\}$, the vertices $x_{i-1 \bmod n}$ and $x_{i+1 \bmod n}$ are neighbours of $x_i$.

Suppose $x_i, x_j \in V_2$ such that $x_ix_j \in E(G)$. For each $\ell \in \{-1,0,1,2\}$, let $y_\ell = x_{i + \ell \bmod n}$. Since $n \geq 4$, the vertices $y_{-1}, y_0, y_1, y_2$ are distinct. Since $y_0 = x_i \in V_2$, we have $N_G(y_0) = \{y_{-1}, y_{1}\}$, so $x_j$ is $y_{-1}$ or $y_{1}$. We may assume that $x_j = y_{1}$. Since $x_j \in V_2$, we obtain $N_G(y_{1}) = \{y_0, y_{2}\}$. Together with $N_G(y_0) = \{y_{-1}, y_{1}\}$, this gives us that the path $y_{-1} y_0 y_1 y_2$ lies on an interior face of $G$, which contradicts the assumption that $G$ is a mop.

Therefore, $V_2$ is an independent set of $G$. Thus, each edge of $G$ contains at most one vertex in $V_2$. Let $H$ be the set of Hamiltonian edges. For any $v \in V_2$ and $e \in H$, let $\chi(v,e) = 1$ if $e$ contains $v$, and let $\chi(v,e) = 0$ otherwise. By the above, $\sum_{v \in V_2} \chi(v,e) \leq 1$ for each $e \in H$. If $x_i \in V_2$, then the edges containing $x_i$ are $x_{i-1 \bmod n}x_i$ and $x_ix_{i+1 \bmod n}$, so $\sum_{e \in H} \chi(x_i,e) = 2$. We have
\[2|V_2| = \sum_{v \in V_2} 2 = \sum_{v \in V_2} \sum_{e \in H} \chi(v,e) = \sum_{e \in H} \sum_{v \in V_2} \chi(v,e) \leq \sum_{e \in H} 1 = n,\]
and hence $|V_2| \leq n/2$.~\hfill{\qed}

The next lemma lists other facts about vertices of degree $2$, which we shall also use.

\begin{lem}\label{basicfacts} If $G$ is a mop of order $n \geq 3$, then the following hold:\\
(a) Each vertex of $G$ is of degree at least $2$. \\
(b) $G$ has at least $2$ vertices of degree $2$. \\
(c) $G - v$ is a mop for any vertex $v$ of $G$ of degree $2$. \\
(d) A graph $H$ is a mop if $G = H-w$ for some $w \in V(H)$ such that $d_H(w) = 2$ and $N_H(w)$ is a Hamiltonian edge of $G$.
\end{lem}

\proof (a) This is immediate from the fact that $G$ has a Hamiltonian cycle.

\noindent
(b) We use induction on $n$. The result is trivial if $n = 3$. Suppose $n \geq 4$. Since $G$ is a mop, $G$ has a diagonal $d = xy$. By Lemma~\ref{diagonallemma}, $G$ is the union of two mops $G_1$ and $G_2$ such that $V(G_1) \cap V(G_2) = d$ and $E(G_1) \cap E(G_2) = \{d\}$. By the induction hypothesis, for each $i \in \{1, 2\}$, $G_i$ has two vertices $v_{i,1}$ and $v_{i,2}$ such that $d_{G_i}(v_{i,1}) = d_{G_i}(v_{i,2}) = 2$. By Lemma~\ref{maxofn2}, for each $i \in \{1, 2\}$, we cannot have $d_{G_i}(x) = d_{G_i}(y) = 2$, so $v_{i,j_i} \notin \{x, y\}$ for some $j_i \in \{1, 2\}$. Therefore, we have $v_{1,j_1} \neq v_{2,j_2}$, $d_G(v_{1,j_1}) = d_{G_1}(v_{1,j_1}) = 2$ and $d_G(v_{2,j_2}) = d_{G_2}(v_{2,j_2}) = 2$.

\noindent
(c) Let $v$ be a vertex of $G$ of degree $2$. We may label the vertices $x_1, x_2, \dots, x_n$ so that $x_1x_2 \dots x_nx_1$ is the Hamiltonian cycle $C$ of $G$ and $x_n = v$. Since $d_G(x_n) = 2$, $N_G(x_n) = \{x_1, x_{n-1}\}$. The face $F$ having $x_{n-1}x_n$ and $x_nx_1$ on its boundary must also have $x_{n-1}x_1$ on its boundary (as all interior faces are triangles), meaning that $x_{n-1}x_1 \in E(G)$. Thus, $x_1x_2 \dots, x_{n-1}x_1$ is a Hamiltonian cycle of $G-v$. Also, every interior face of $G$ other than $F$ is a face of $G-v$ (and a triangle). Therefore, $G-v$ is a mop.

\noindent
(d) Let $x_1x_2 \dots x_nx_1$ be the Hamiltonian cycle $C$ of $G$. We may assume that $N_H(w) = \{x_1, x_2\}$. Let $C' = x_1wx_2 \dots x_nx_1$. Then, $C'$ is a Hamiltonian cycle of $H$. Let $\rho$ be a plane drawing of $G$ such that all vertices of $G$ lie on the boundary $C$ and all interior faces are triangles. Extend $\rho$ to a plane drawing $\rho'$ of $H$ by putting $w$ and the edges $wx_1, wx_2$ on the exterior of $\rho$. Then, $C'$ is the boundary of $\rho'$. Also, the faces of $\rho'$ are the faces of $\rho$ together with the face bounded by $wx_1$, $wx_2$ and $x_1x_2$, so all faces of $\rho'$ are triangles. Thus, $H$ is a mop.\hfill{\qed}

Lemmas~\ref{maxofn2} and \ref{basicfacts}(b) tell us that the number $n_2$ of vertices of degree $2$ (of a mop) satisfies
\begin{equation} 2 \leq n_2 \leq \frac{n}{2}. \label{n_2bounds}
\end{equation}
We show in passing that both bounds are sharp. For the upper bound, we start with any mop $M$ with Hamiltonian cycle $x_1x_2 \dots x_px_1$, add $p$ new vertices $y_1, y_2, \dots, y_p$, and add the edges $y_1x_1, y_1x_2$, $y_2x_2, y_2x_3$, $\dots, y_px_p, y_px_1$, and hence the resulting graph $G$ is a mop (by Lemma~\ref{basicfacts}(d)), the vertices of $G$ of degree $2$ are $y_1, y_2, \dots, y_p$, and $p = \frac{|V(G)|}{2}$. For the lower bound, we start with a cycle $C = x_1x_2 \dots x_nx_1$, set $p = \lceil n/2 \rceil$, and take $G$ with $V(G) = V(C)$ and $E(G) = E(C)\cup \{x_2x_n, x_nx_3, x_3x_{n-1}, x_{n-1}x_4, \dots, x_{p-1}x_{n-p+3}, x_{n-p+3}x_p, x_px_{n-p+2}\}$, and hence the vertices of $G$ of degree $2$ are $x_1$ and $x_{p+1}$.

The next lemma settles Theorem~\ref{S:main1} for $5 \leq n \leq 9$, and hence allows us to use Lemma~\ref{mop:5-8} in the proof of Theorem~\ref{S:main1}.

\begin{lem}\label{natmost9}
If $5 \leq n \le 9$ and $G$ is a mop of order $n$, then $\ii_1(G) = 1 \leq \frac{n}{5}$.
\end{lem}
\proof
By definition of a mop, $\ii_1(G) \ge 1$.

Suppose $n = 9$. Let $x_1x_2 \dots x_9x_1$ be the unique Hamiltonian cycle of $G$ and hence the boundary of the exterior face of $G$. By Lemma \ref{mop:4-6}, $G$ has a diagonal $d$ such that $G$ is the union of two mops $G_1$ and $G_2$ such that $V(G_1) \cap V(G_2) = d$, $E(G_1) \cap E(G_2) = \{d\}$, and $G_1$ has exactly $\ell$ Hamiltonian edges of $G$ for some $\ell \in \{4,5,6\}$. We may assume that $d$ is the edge $x_1x_{\ell+1}$ and that $V(G_1) = \{x_1, x_2, \dots, x_{\ell+1}\}$.

Suppose $\ell = 4$. Then, $V(G_{2})= \{x_1, x_5,x_6,x_7,x_8,x_9\}$. Let $(x_1,x_5,x_j)$ be the triangular face of $G_2$ containing the edge $x_1x_5$ (so $j \in \{6,7,8,9\}$). Let $x'=x_5$ if $j=9$, and let $x'=x_1$ otherwise. If $6 \le j \le 8$, then $x_1,x_2,x_5,x_j,x_9 \in N_{G}[x']$. If $j=9$, then $x_1,x_4,x_5,x_6,x_9 \in N_{G}[x']$. Clearly, each component of $G - N_{G}[x']$ contains at most two vertices. Therefore, $\{x'\}$ is a $K_{1,2}$-isolating set of $G$. Hence, $\ii_1(G) = 1 < n/5$.

Suppose $\ell = 5$. Then, $x_6x_7,x_7x_8,x_8x_9,x_9x_1$ are the $4$ Hamiltonian edges of $G$ that belong to $G_2$. Thus, this case is identical to the case $\ell = 4$.

Suppose $\ell = 6$. Let $(x_1,x_j,x_7)$ be the triangular face of $G_1$ containing the edge $x_1x_7$ (so $j \in \{2,3,4,5,6\}$). Suppose $j=2$. Then, $x_2x_7 \in E(G)$. By Lemma~\ref{diagonallemma}, $x_2x_7$ partitions $G$ into two mops $H_1$ and $H_2$ such that $H_1$ has exactly $4$ Hamiltonian edges of $G$ (namely, $x_7x_8, x_8x_9, x_9x_1, x_1x_2$), so the result follows as in the case $\ell = 4$. By symmetry, it suffices to consider the cases $j = 3$ and $j = 4$; that is, the cases $j = 6$ and $j = 5$ are identical to the cases $j = 2$ and $j = 3$, respectively. In both cases, each component of $G - N_{G}[x_j]$ contains at most two vertices, so $\{x_j\}$ is a $K_{1,2}$-isolating set of $G$.

Now suppose $n = 8$. Let $uv$ be a Hamiltonian edge of $G$, and let $H$ be the graph obtained by adding a new vertex $w$ to $G$ and adding the edges $wu$ and $wv$. By Lemma~\ref{basicfacts}(d), $H$ is a mop of order $9$, so $\ii_1(H) = 1 < n/5$. Let $\{x\}$ be a $K_{1,2}$-isolating set of $H$. If $x \neq w$, then $\{x\}$ is a $K_{1,2}$-isolating set of $G$. If $x = w$, then $\{u\}$ is a $K_{1,2}$-isolating set of $G$ as $N_{H}[w] = \{u,v,w\} \subseteq N_{H}[u]$. Therefore, $\ii_1(G) = 1 < n/5$.

For $5 \leq i \leq 7$, we obtain the result for $n = i$ from the result for $n = i+1$ in the same way we obtained the result for $n = 8$ from the result for $n = 9$.
\hfill{\qed}

We now prove Theorems~\ref{t:main1} and~\ref{t:main2}. Recall the statement of Theorem~\ref{t:main1}.

\noindent \textbf{Theorem~\ref{t:main1}}. \emph{If $G$ is a mop of order $n \ge 5$, then $\ii_{1}(G) \leq \frac{n}{5}$.}

\proof We use induction on $n$. Lemma~\ref{natmost9} establishes the base case $n = 5$ (in this case, $G$ is a fan $F_5$, and its center is a $K_{1,2}$-isolating set of $G$) and also the case $6 \leq n \leq 9$. Suppose $n \ge 10$. We assume that if $G'$ is a mop of order $n'$ with $5 \le n' < n$, then $\ii_1(G') \le n'/5$.

Let $x_1x_2 \ldots x_n x_1$ be the unique Hamiltonian cycle $C$ of $G$ and hence the boundary of the exterior face of $G$. By Lemma~\ref{mop:5-8}, $G$ has a diagonal $d$ such that $G$ is the union of two mops $G_1$ and $G_2$ such that $V(G_1) \cap V(G_2) = d$, $E(G_1) \cap E(G_2) = \{d\}$, and $G_1$ has exactly $\ell$ Hamiltonian edges of $G$ for some $\ell \in \{5,6,7,8\}$. We may assume that $d$ is the edge $x_1x_{\ell + 1}$ and that $V(G_1) =  \{x_1,x_2,\ldots,x_{\ell + 1}\}$. Note that the edges $x_1x_2,x_2x_3,\dots,x_{\ell}x_{\ell + 1}$ are the $\ell$ Hamiltonian edges of $G$ that belong to $G_1$. Let $(x_1,x_j,x_{\ell + 1})$ be the triangular face of $G_1$ containing the edge $x_1x_{\ell + 1}$ (so $j \in \{2,3,\dots,\ell\}$).

\begin{claim}
\label{c1:1}
If $\ell = 5$, then $\ii_1(G) \le \frac{n}{5}$.
\end{claim}
\proof Let $G'$ be the graph obtained from $G$ by deleting the vertices $x_2,x_3,x_4,x_5$ and contracting the edge $x_1x_6$ to form a new vertex $y$ (see Figure 4). Thus, $G'$ is obtained from $G_2$ by contracting the edge $x_1x_6$. By Lemma~\ref{contraction}, $G'$ is a mop.
\vskip 5 pt

\setlength{\unitlength}{0.6cm}
\begin{center}
\begin{picture}(15, 4)
\put(1, 0){\circle*{0.2}}
\put(3, 0){\circle*{0.2}}
\put(3, 3){\circle*{0.2}}
\put(1, 3){\circle*{0.2}}
\put(4.5, 1.5){\circle*{0.2}}
\put(-0.5, 1.5){\circle*{0.2}}

\put(-1.2, 1.3){\small{$x_2$}}
\put(1, 3.2){\small{$x_3$}}
\put(3, 3.2){\small{$x_4$}}
\put(4.7, 1.3){\small{$x_5$}}
\put(-0.2, 0.2){\small{$x_1$}}
\put(3.5, 0.2){\small{$x_6$}}
\put(-2, 0.2){\small{$x_n$}}
\put(6, 0.2){\small{$x_7$}}
\put(-2, -3.5){\small{$x_p$}}
\put(2, -3.5){\small{$x_q$}}
\put(6, -3.5){\small{$x_r$}}

\put(1, 0){\line(-1, 1){1.5}}
\put(3, 0){\line(1, 1){1.5}}
\put(4.5, 1.5){\line(-1, 1){1.5}}
\put(-0.5,1.5){\line(1, 1){1.5}}
\put(1, 3){\line(1, 0){2}}
\put(1, 0){\line(1, 0){2}}
\put(-2, 0){\circle*{0.2}}
\put(6, 0){\circle*{0.2}}
\put(2, -3){\circle*{0.2}}
\put(-2, -3){\circle*{0.2}}
\put(6, -3){\circle*{0.2}}
\put(6, -3){\line(-1, 1){3}}
\put(-2, -3){\line(1, 1){3}}
\put(2, -3){\line(-1, 3){1}}
\put(2, -3){\line(1, 3){1}}
\put(-2, 0){\line(1, 0){8}}
\put(0, -2){\small{$\dots$}}
\put(3, -2){\small{$\dots$}}
\put(10, 0.2){\small{$x_n$}}
\put(18, 0.2){\small{$x_6$}}
\put(14, 0.2){\small{$y$}}
\put(10, -3.5){\small{$x_p$}}
\put(18, -3.5){\small{$x_r$}}
\put(14, -3.5){\small{$x_q$}}
\put(14, 0){\circle*{0.2}}
\put(10, 0){\circle*{0.2}}
\put(18, 0){\circle*{0.2}}
\put(14, -3){\circle*{0.2}}
\put(10, -3){\circle*{0.2}}
\put(18, -3){\circle*{0.2}}
\put(18, -3){\line(-4, 3){4}}
\put(10, -3){\line(4, 3){4}}
\put(14, -3){\line(0, 1){3}}
\put(10, 0){\line(1, 0){8}}
\put(12, -2){\small{$\dots$}}
\put(15, -2){\small{$\dots$}}
\put(7.5, 0){\Large{$\Rightarrow$}}

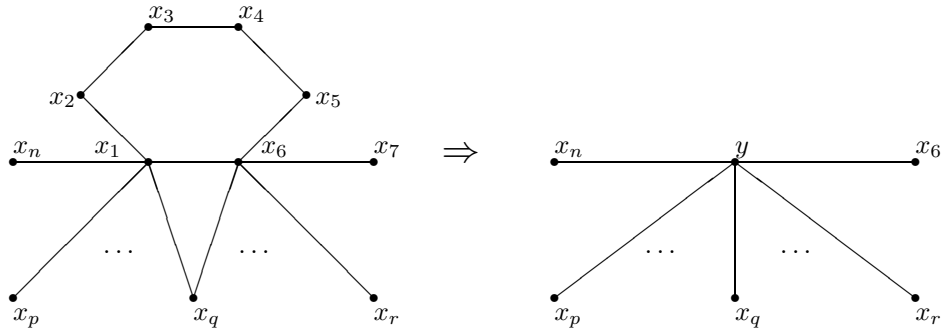
\put(3.2, -5){\footnotesize{\textbf{\textbf{Figure 4 :}} The edge contraction of $G_2$.}}
\end{picture}
\end{center}
\vskip 90 pt

Let $n'$ be the order $n - 5 \geq 5$ of $G'$. By the induction hypothesis, $\ii_1(G') \le n'/5 = n/5-1$. Let $S'$ be a smallest $K_{1,2}$-isolating set of $G'$. Then, $|S'| = \ii_1(G') \le n/5-1$. We have that $\{x_1,x_6\}$ is a $\{K_{1,2}\}$-isolating set of $G_1$ and that $N_{G'}(y) \subseteq N_{G}(x_1) \cup N_{G}(x_6)$. Thus, if $y \in S'$, then $(S' \setminus \{y\}) \cup \{x_1,x_6\}$ is a $\{K_{1,2}\}$-isolating set of $G$, and hence $\ii_1(G) \le (|S'| - 1) + 2 \le n/5$. Suppose $y \notin S'$. Observe that $|V(G_1 - N_{G_1}[x_j])| \le 2$. Thus, $G_1 - N_{G_1}[x_j]$ does not contain a copy of $K_{1,2}$. Moreover, $S' \cup \{x_j\}$ is a $K_{1,2}$-isolating set of $G$ as $x_j$ is adjacent to both $x_1$ and $x_6$. Therefore, $\ii_1(G) \le |S'| + 1 \le n/5$.~\hfill{\smallqed}

\begin{claim}
\label{c1:2}
If $\ell = 6$, then $\ii_1(G) \le \frac{n}{5}$.
\end{claim}
\proof Let $G' = G_2$ and $n' = |V(G_2)|$. We have $V(G') = V(G) \setminus \{x_2,x_3,x_4,x_5,x_6\}$, so $n' = n-5 \geq 5$. By the induction hypothesis, $\ii_1(G') \le n'/5 = n/5-1$. Let $S'$ be a smallest $K_{1,2}$-isolating set of $G'$. Then, $|S'| = \ii_1(G') \le n/5-1$. Suppose $j=2$. Then, $x_2x_7$ is a diagonal of $G$ that partitions $G$ into two mops $H_1$ and $H_2$ with $H_1 = G_1 - x_1$. Since $H_1$ has exactly $5$ Hamiltonian edges of $G$ (namely, $x_2x_3,x_3x_4,x_4x_5,x_5x_6,x_6x_7$), we have $\ii_1(G) \le \frac{n}{5}$ by Claim \ref{c1:1}. Now suppose $j \ge 3$. By symmetry, it suffices to consider the cases $j = 3$ and $j = 4$ (see Figure 5); that is, the cases $j = 6$ and $j = 5$ are identical to the cases $j = 2$ and $j = 3$, respectively. In both cases, the order of each component of $G_1 - N_{G_1}[x_j]$ is at most $2$, so $G_1 - N_{G_1}[x_j]$ does not contain a copy of $K_{1,2}$. Since $x_j$ is adjacent to $x_1$ and $x_7$, $S' \cup \{x_j\}$ is a $K_{1,2}$-isolating set of $G$, and hence $\ii_1(G) \le |S'| + 1 \le n/5$.~\hfill{\smallqed}

\setlength{\unitlength}{0.6cm}
\begin{center}
\begin{picture}(17, 6)
\put(1, 0){\circle*{0.2}}
\put(5, 0){\circle*{0.2}}
\put(0, 2){\circle*{0.2}}
\put(1, 4){\circle*{0.2}}
\put(5, 4){\circle*{0.2}}
\put(6, 2){\circle*{0.2}}
\put(3, 5){\circle*{0.2}}
\put(1, 0){\line(0, 1){4}}
\put(-1, 0){\line(1, 0){8}}
\put(5, 0){\line(-1, 1){4}}
\put(1, 0){\line(-1, 2){1}}
\put(1, 4){\line(-1, -2){1}}
\put(1, 4){\line(2, 1){2}}
\put(5, 4){\line(-2, 1){2}}
\put(6, 2){\line(-1, 2){1}}
\put(6, 2){\line(-1, -2){1}}
\put(0, 0.2){\small{$x_1$}}
\put(-0.8, 2){\small{$x_2$}}
\put(0.6, 4.2){\small{$x_3$}}
\put(3, 5.2){\small{$x_4$}}
\put(5, 4.2){\small{$x_5$}}
\put(6.2, 2){\small{$x_6$}}
\put(5.4, 0.2){\small{$x_7$}}
\put(13, 0){\circle*{0.2}}
\put(17, 0){\circle*{0.2}}
\put(12, 2){\circle*{0.2}}
\put(13, 4){\circle*{0.2}}
\put(17, 4){\circle*{0.2}}
\put(18, 2){\circle*{0.2}}
\put(15, 5){\circle*{0.2}}
\put(15, 5){\line(2, -5){2}}
\put(15, 5){\line(-2, -5){2}}
\put(11, 0){\line(1, 0){8}}
\put(13, 0){\line(-1, 2){1}}
\put(13, 4){\line(-1, -2){1}}
\put(13, 4){\line(2, 1){2}}
\put(17, 4){\line(-2, 1){2}}
\put(18, 2){\line(-1, 2){1}}
\put(18, 2){\line(-1, -2){1}}
\put(12, 0.2){\small{$x_1$}}
\put(11.2, 2){\small{$x_2$}}
\put(12.6, 4.2){\small{$x_3$}}
\put(15, 5.2){\small{$x_4$}}
\put(17, 4.2){\small{$x_5$}}
\put(18.2, 2){\small{$x_6$}}
\put(17.4, 0.2){\small{$x_7$}}

\put(7.5, -2){\footnotesize{\textbf{\textbf{Figure 5  }}  }}
\end{picture}
\end{center}
\vskip 50 pt

\begin{claim}
\label{c1:3}
If $\ell = 7$, then $\ii_1(G) \le \frac{n}{5}$.
\end{claim}
\proof Let $G' = G_2$ and $n' = |V(G_2)|$. We have $V(G') = V(G) \setminus \{x_2,x_3,x_4,x_5,x_6,x_7\}$, so $n' = n-6 \ge 4$. Suppose $j=2$. Then, $x_2x_8$ is a diagonal of $G$ that partitions $G$ into two mops $H_1$ and $H_2$ with $H_1 = G_1 - x_1$. Since $H_1$ has exactly $6$ Hamiltonian edges of $G$ (namely, $x_2x_3,x_3x_4,\dots,x_7x_8$), we have $\ii_1(G) \le \frac{n}{5}$ by Claim \ref{c1:2}. Similarly, if $j=3$, then $x_3x_8$ is a diagonal which partitions $G$ into two mops $H_1$ and $H_2$ with $H_1 = G_1 - \{x_1, x_2\}$, and hence the result follows by Claim \ref{c1:1}. Now suppose $j \ge 4$. By symmetry, we may assume that $j = 4$ (see Figure 6, left); that is, the cases $j = 7$, $j = 6$ and $j = 5$ are identical to the cases $j = 2$, $j = 3$ and $j = 4$, respectively.

The order of each component of $G_1 - N_{G_1}[x_4]$ is at most $2$, so $G_1 - N_{G_1}[x_4]$ does not contain a copy of $K_{1,2}$. Suppose $n' = 4$. Then, $n = 10$ and $|V(G')| = 4$. Since $x_4$ is adjacent to both $x_1$ and $x_8$, $G' - N_{G}[x_4]$ does not contain a copy of $K_{1,2}$. Thus, $\{x_4\}$ is a $K_{1,2}$-isolating set of $G$, and hence $\ii_1(G) = 1 < n/5$. Now suppose $n' \ge 5$. By the induction hypothesis, $\ii_1(G') \le n'/5 < n/5-1$. Let $S'$ be a smallest $K_{1,2}$-isolating set of $G'$. Then, since $S' \cup \{x_4\}$ is a $K_{1,2}$-isolating set of $G$, we have $\ii_1(G) \le |S'| + 1 = \ii_1(G') + 1 < n/5$.~\hfill{\smallqed}

\setlength{\unitlength}{0.6cm}
\begin{center}
\begin{picture}(20, 7)
\put(0, 0){\circle*{0.2}}
\put(6, 0){\circle*{0.2}}
\put(-1, 2){\circle*{0.2}}
\put(0, 4){\circle*{0.2}}
\put(1.8, 5.35){\circle*{0.2}}
\put(4.2, 5.35){\circle*{0.2}}
\put(6, 4){\circle*{0.2}}
\put(7, 2){\circle*{0.2}}
\put(0, 0){\line(1, 3){1.8}}
\put(6, 0){\line(-4, 5){4.3}}
\put(-2, 0){\line(1, 0){10}}
\put(0, 0){\line(-1, 2){1}}
\put(0, 4){\line(-1, -2){1}}
\put(0, 4){\line(4, 3){1.8}}
\put(6, 4){\line(-4, 3){1.8}}
\put(7, 2){\line(-1, 2){1}}
\put(7, 2){\line(-1, -2){1}}
\put(1.8, 5.35){\line(1, 0){2.3}}
\put(-1.2, 0.2){\small{$x_1$}}
\put(-2, 2){\small{$x_2$}}
\put(-0.6, 4.2){\small{$x_3$}}
\put(1.6, 5.5){\small{$x_4$}}
\put(4.4, 5.5){\small{$x_5$}}
\put(6.2, 4.2){\small{$x_6$}}
\put(7.4, 2){\small{$x_7$}}
\put(6.6, 0.2){\small{$x_8$}}

\put(13, 0){\circle*{0.2}}
\put(19, 0){\circle*{0.2}}
\put(12, 2){\circle*{0.2}}
\put(13, 4){\circle*{0.2}}
\put(14.8, 5.3){\circle*{0.2}}
\put(17.2, 5.3){\circle*{0.2}}
\put(19, 4){\circle*{0.2}}
\put(20, 2){\circle*{0.2}}
\put(16, 6.15){\circle*{0.2}}

\put(13, 0.1){\line(1, 2){3}}
\put(19, 0.1){\line(-1, 2){3}}

\put(11, 0){\line(1, 0){10}}
\put(13, 0){\line(-1, 2){1}}
\put(13, 4){\line(-1, -2){1}}
\put(13, 3.95){\line(4, 3){3}}
\put(19, 3.95){\line(-4, 3){3}}
\put(20, 2){\line(-1, 2){1}}
\put(20, 2){\line(-1, -2){1}}

\put(11.8, 0.2){\small{$x_1$}}
\put(11, 2){\small{$x_2$}}
\put(12.4, 4.2){\small{$x_3$}}
\put(14.2, 5.5){\small{$x_4$}}
\put(16, 6.5){\small{$x_5$}}
\put(17.4, 5.5){\small{$x_6$}}
\put(19.2, 4.2){\small{$x_7$}}
\put(20.4, 2){\small{$x_8$}}
\put(19.6, 0.2){\small{$x_9$}}

\put(8.5, -2){\footnotesize{\textbf{\textbf{Figure 6 }} }}
\end{picture}
\end{center}
\vskip 50 pt

By Claims~\ref{c1:1},~\ref{c1:2} and ~\ref{c1:3}, we may assume that $\ell = 8$. Similarly to the above, the cases $j = 2$, $j = 3$ and $j = 4$ follow by Claims \ref{c1:3}, \ref{c1:2} and \ref{c1:1}, respectively. Thus, we may assume that $j \ge 5$. Also similarly to the above, by symmetry, it suffices to consider $j = 5$ (see Figure 6, right).

The order of each component of $G_1 - N_{G_1}[x_5]$ is at most $2$, so $G_1 - N_{G_1}[x_5]$ does not contain a copy of $K_{1,2}$. Let $G' = G_2$ and $n' = |V(G_2)|$. Then, $n' = n - 7$.  Suppose $n' \le 4$. Then, $n \le 11$ and $|V(G')| \le 4$. Since $x_5$ is adjacent to both $x_1$ and $x_9$, $G' - N_{G}[x_5]$ does not contain a copy of $K_{1,2}$. Thus, $\{x_5\}$ is a $K_{1,2}$-isolating set of $G$, and hence $\ii_1(G) = 1 < n/5$. Now suppose $n' \ge 5$. By the induction hypothesis, $\ii_1(G') \le n'/5 < n/5-1$. Let $S'$ be a smallest $K_{1,2}$-isolating set of $G'$. Then, since $S' \cup \{x_5\}$ is a $K_{1,2}$-isolating set of $G$, we have $\ii_1(G) \le |S'| + 1 = \ii_1(G') + 1 < n/5$.~\hfill{\qed}

Recall the statement of Theorem~\ref{t:main2}.

\noindent \textbf{Theorem~\ref{t:main2}}. \emph{If $G$ is a mop of order $n \ge 5$ with $n_2$ vertices of degree~$2$, then
\[ \ii_1(G)  \le \left\{\begin{array}{ll}
                     \frac{n + n_2}{6} & when~~n_2 \le \frac{n}{3},\\
                     \\
                     \frac{n - n_2}{3} & otherwise.
               \end{array}\right. \]
}

\proof We first consider $n_2 > \frac{n}{3}$. Let $V_2$ be the set of vertices of $G$ of degree~$2$. Let $G' = G - V_2$ and $n' = |V(G')|$. By Lemma~\ref{maxofn2}, $V_2$ is an independent set of $G$, and $n_2 \leq n/2$. Since $n'$ is an integer satisfying $n' = n - n_2 \geq n/2 \geq 5/2$, we have $n' \geq 3$. Since $G'$ is obtained from $G$ by deleting the vertices in $V_2$ one by one, $G'$ is a mop by Lemma~\ref{basicfacts}(c). Let $S$ be a smallest dominating set of $G'$. By Theorem~\ref{t:known1}(a), $|S| \le n'/3$. Since $V_2$ is an independent set of $G$, $S$ is a $K_{1, 2}$-isolating set of $G$. Thus, $\ii_1(G) \le |S| \le (n - n_2)/3$.

We now prove the bound for $n_2 \le \frac{n}{3}$, using induction on $n$. We actually prove it for $n \ge 4$. If $n = 4$, then we trivially have $\ii_1(G) = 1$, so $\ii_1(G) \le (n+n_2)/6$ by Lemma~\ref{basicfacts}(b). By Theorem~\ref{t:main1} and Lemma~\ref{basicfacts}(b), we have $\ii_1(G) \le 1 \le (n+n_2)/6$ if $5 \leq n \le 9$, and we have $\ii_1(G) \le 2 \le (n+n_2)/6$ if $n = 10$. We now consider $n \ge 11$ and assume that $\ii_1(G') \le (n'+n_2')/6$ for any mop $G'$ such that $|V(G')| = n' \geq 4$, $G'$ has $n_2'$ vertices of degree~$2$, and $n'+n_2' < n+n_2$.

Let $x_1x_2 \ldots x_n x_1$ be the unique Hamiltonian cycle $C$ of $G$ and hence the boundary of the exterior face of $G$. By Lemma~\ref{mop:5-8}, $G$ has a diagonal $d$ such that $G$ is the union of two mops $G_1$ and $G_2$ such that $V(G_1) \cap V(G_2) = d$, $E(G_1) \cap E(G_2) = \{d\}$, and $G_1$ has exactly $\ell$ Hamiltonian edges of $G$ for some $\ell \in \{5,6,7,8\}$. We may assume that $d$ is the edge $x_1x_{\ell + 1}$ and that $V(G_1) =  \{x_1,x_2,\ldots,x_{\ell + 1}\}$. Note that the edges $x_1x_2,x_2x_3,\dots,x_{\ell}x_{\ell + 1}$ are the $\ell$ Hamiltonian edges of $G$ that belong to $G_1$. Let $(x_1,x_j,x_{\ell + 1})$ be the triangular face of $G_1$ containing the edge $x_1x_{\ell + 1}$ (so $j \in \{2,3,\dots,\ell\}$).

\begin{claim}
\label{c2:1}
If $\ell = 5$, then $\ii_1(G) \le \frac{n+n_2}{6}$.
\end{claim}
\proof By Lemma~\ref{basicfacts}(a), $d_{G_2}(v) \geq 2$ for each $v \in V(G_2)$. Since $x_1$ and $x_6$ are adjacent in $G_2$, Lemma~\ref{maxofn2} tells us that their degrees cannot be both $2$. We may assume that $d_{G_2}(x_1) \geq 3$. This yields $d_{G_2}(x_1) + d_{G_2}(x_6) \ge 5$.

\setlength{\unitlength}{0.8cm}
\begin{center}
\begin{picture}(4, 4)
\put(1, -2){\circle*{0.18}}
\put(3, -2){\circle*{0.18}}
\put(1, 0){\circle*{0.18}}
\put(3, 0){\circle*{0.18}}
\put(1, 3){\circle*{0.18}}
\put(3, 3){\circle*{0.18}}
\put(3, 0){\circle*{0.18}}
\put(4.5, 1.5){\circle*{0.18}}
\put(-0.5, 1.5){\circle*{0.18}}

\put(0.3, 0){\small{$x_1$}}
\put(3.3, 0){\small{$x_6$}}
\put(-1.2, 1.3){\small{$x_2$}}
\put(1, 3.2){\small{$x_3$}}
\put(3, 3.2){\small{$x_4$}}
\put(4.7, 1.3){\small{$x_5$}}
\put(0.3, -1.7){\small{$x_n$}}
\put(3.3, -1.7){\small{$x_7$}}

\put(1, 0){\line(1, 0){2}}
\put(1, 0){\line(-1, 1){1.5}}
\put(3, 0){\line(1, 1){1.5}}
\put(4.5, 1.5){\line(-1, 1){1.5}}
\put(-0.5,1.5){\line(1, 1){1.5}}
\put(1, 3){\line(1, 0){2}}
\put(1, 0){\line(1, 0){2}}
\put(1, -2){\line(0, 1){2}}
\put(-2, -2){\line(1, 0){8}}
\put(3, -2){\line(0, 1){2}}
\put(3, -2){\line(-1, 1){2}}
\put(1.25, -3.3){\footnotesize{\textbf{\textbf{Figure 7}}}}
\end{picture}
\end{center}
\vskip 90 pt

\begin{subclaim}
\label{c2:11}
If $d_{G_2}(x_1) + d_{G_2}(x_6) = 5$, then $\ii_1(G) \le \frac{n+n_2}{6}$.
\end{subclaim}
\proof Suppose $d_{G_2}(x_1) + d_{G_2}(x_6) = 5$. Then, $d_{G_2}(x_1) = 3$ and $d_{G_2}(x_6) = 2$. We have $x_1x_n, x_6x_7 \in E(C) \cap E(G_2)$. Since $x_1x_6, x_6x_7 \in E(G_2)$ and $d_{G_2}(x_6) = 2$, $N_{G_2}(x_6) = \{x_1,x_7\}$. Thus, since $G_2$ is a mop, the face having $x_1x_6$ and $x_6x_7$ on its boundary must also have $x_1x_7$ on its boundary (as all interior faces are triangles), that is, $x_1x_7 \in E(G_2)$ (see Figure 7). Together with $x_1x_6, x_1x_n \in E(G_2)$ and $d_{G_2}(x_1) = 3$, this gives us $N_{G_2}(x_1) = \{x_6,x_7,x_n\}$. Thus, since $G_2$ is a mop, the face having $x_1x_7$ and $x_1x_n$ on its boundary must also have $x_7x_n$ on its boundary, that is, $x_7x_n \in E(G_2)$ (see Figure 7). Let $G' = G - \{x_1,x_2,\dots,x_6\}$. Then, $G' = (G_2 - x_6) - x_1$. Since $d_{G_2}(x_6) = 2$, $G_2 - x_6$ is a mop by Lemma~\ref{basicfacts}(c). Since $d_{G_2 - x_6}(x_1) = 2$, $G'$ is a mop by Lemma~\ref{basicfacts}(c). Let $n' = |V(G')|$ and $n_2' = |\{v \in V(G') \colon d_{G'}(v) = 2\}|$. We have $n' = n - 6 \ge 5$.

By Lemma~\ref{maxofn2}, at most one of $x_7$ and $x_n$ has degree~$2$ in $G'$. By Lemma~\ref{maxofn2}, at most one of $x_1$ and $x_6$ has degree~$2$ in $G_1$, and hence, by Lemma~\ref{basicfacts}(b), $d_{G_1}(x_h) = 2$ for some $h \in \{2, 3, 4, 5\}$. Since $x_h \in V(G_1) \backslash V(G_2)$, $d_{G}(x_h) = d_{G_1}(x_h)$. Therefore, $n_2' \le n_2$, and hence $n'+n_2' \le n+n_2-6$. By the induction hypothesis, $\ii_1(G') \le (n'+n_2')/6 \le (n+n_2)/6-1$. Let $S'$ be a smallest $K_{1,2}$-isolating set of $G'$. Clearly, $|V(G_1 - N_{G_1}[x_j])| \le 2$, so $G_1 - N_{G_1}[x_j]$ does not contain a copy of $K_{1, 2}$. Since $x_j$ is adjacent to both $x_1$ and $x_6$, it follows that $S' \cup \{x_j\}$ is a $K_{1,2}$-isolating set of $G$. Thus, we have $\ii_1(G) \le |S'| + 1 = \ii_1(G') + 1 \le (n+n_2)/6$.~\hfill{\smallqed}
\medskip

By Claim~\ref{c2:11}, we may assume that $d_{G_2}(x_1) + d_{G_2}(x_6) \ge 6$. Let $G'$ be the graph obtained from $G$ by deleting the vertices $x_2,x_3,x_4,x_5$ and contracting the edge $x_1x_6$ to form a new vertex $y$. Then, $G'$ is obtained from $G_2$ by contracting $x_1x_6$. Thus, $G'$ is a mop by Lemma~\ref{contraction}. Let $n' = |V(G')|$ and $n_2' = |\{v \in V(G') \colon d_{G'}(v) = 2\}|$. We have $n' = n - 5 \ge 6$.

Suppose that every vertex that has degree~$2$ in $G'$ also has degree~$2$ in $G$. As in the proof of Claim~\ref{c2:11}, $d_{G_1}(x_h) = 2$ for some $h \in \{2, 3, 4, 5\}$, so $n_2' \le n_2 - 1$. Thus, $n'+n_2' \le n+n_2-6$. By the induction hypothesis, $\ii_1(G') \le (n'+n_2')/6 \le (n+n_2)/6-1$. Let $S'$ be a smallest $K_{1,2}$-isolating set of $G'$. Then, $|S'| \le (n+n_2)/6-1$.
We can continue as in the proof of Claim~\ref{c1:1} to obtain $\ii_1(G) \le |S'| + 1 \le (n+n_2)/6$.

Now suppose that $G'$ has a vertex $z$ such that $d_{G'}(z) = 2 \neq d_G(z)$. Then $z = y$ or $x_1, x_6 \in N_{G_2}(z)$.

Suppose $z = y$. As noted in the proof of Claim~\ref{c2:11}, $x_1x_n, x_6x_7 \in E(G_2)$. Thus, $x_7, x_n \in N_{G'}(y)$. Since $d_{G'}(y) = 2$, $N_{G'}(y) = \{x_7, x_n\}$. Since $d_{G_2}(x_1) + d_{G_2}(x_6) \ge 6$, we obtain $N_{G_2}(x_1) = \{x_6, x_7, x_n\}$ and $N_{G_2}(x_6) = \{x_1, x_7, x_n\}$. Since $N_{G_2}(x_1) = \{x_6, x_7, x_n\}$, $x_1x_7$ is a diagonal of $G_2$. By Lemma~\ref{diagonallemma}, we obtain $x_6x_n \notin E(G_2)$, which contradicts $N_{G_2}(x_6) = \{x_1, x_7, x_n\}$.

Therefore, $z \neq y$. Thus, $x_1, x_6 \in N_{G_2}(z)$. For each $i \in \{8, 9, \dots, n-1\}$ with $d_{G'}(x_i) = 2$, we have $N_{G'}(x_i) = N_G(x_i) = \{x_{i-1}, x_{i+1}\}$, so $z \neq x_i$. Thus, $z = x_7$ or $z = x_n$. By symmetry, we may assume that $z = x_7$. Since $x_7x_8 \in E(C) \cap E(G_2)$ and $x_1, x_6 \in N_{G_2}(x_7)$, $N_{G_2}(x_7) = \{x_1, x_6, x_8\}$. Thus, since $G_2$ is a mop, the face having $x_1x_7$ and $x_7x_8$ on its boundary must also have $x_1x_8$ on its boundary (as all interior faces are triangles), meaning that $x_1x_8 \in E(G_2)$. By Lemma~\ref{diagonallemma}, $x_1x_8$ partitions $G$ into two mops $H_1$ and $H_2$ such that $V(H_2) = \{x_1, x_8, x_9, \dots, x_n\}$. Let $G' = H_2$, $n' = V(H_2)$ and $n_2' = |\{v \in V(H_2) \colon d_{H_2}(v) = 2\}|$. We have $n' = n-6 \geq 5$. By Lemma~\ref{maxofn2}, for each $i \in \{1, 2\}$, at most one of $x_1$ and $x_8$ has degree~$2$ in $H_i$. By Lemma~\ref{basicfacts}(b), $d_{H_1}(x_h) = 2$ for some $h \in V(H_1) \backslash \{x_1, x_8\}$, and hence $d_{G}(x_h) = 2$. Therefore, $n_2' \le n_2$, and hence $n'+n_2' \le n+n_2-6$. By the induction hypothesis, $\ii_1(G') \le (n'+n_2')/6 \le (n+n_2)/6-1$. Let $S'$ be a smallest $K_{1,2}$-isolating set of $G'$. Let $x' = x_6$ if $j = 2$, and let $x' = x_1$ otherwise. If $3 \leq j \leq 5$, then $x_1, x_2, x_j, x_6, x_7, x_8 \in N_{H_1}[x']$. If $j = 2$, then $x_1, x_2, x_5, x_6, x_7 \in N_{H_1}[x']$. Since $x_1x_6$ is a diagonal of $G$, we have $x_3, x_4 \notin N_G(x_8)$ by Lemma~\ref{diagonallemma}. Therefore, $S' \cup \{x'\}$ is a $K_{1,2}$-isolating set of $G$. Thus, we have $\ii_1(G) \le |S'| + 1 = \ii_1(G') + 1 \le (n+n_2)/6$.~\hfill{\smallqed}

\begin{claim}
\label{c2:2}
If $\ell = 6$, then $\ii_1(G) \le \frac{n+n_2}{6}$.
\end{claim}
\proof As in the proof of Claim~\ref{c1:2}, we may assume that $j = 3$ or $j = 4$. In both cases, the order of each component of $G_1 - N_{G_1}[x_j]$ is at most $2$, so $G_1 - N_{G_1}[x_j]$ does not contain a copy of $K_{1,2}$. Let $G' = G_2$, $n' = |V(G_2)|$ and $n_2' = |\{v \in V(G_2) \colon d_{G_2}(v) = 2\}|$. We have $V(G') = V(G) \setminus \{x_2,x_3,x_4,x_5,x_6\}$, so $n' = n-5 \geq 6$.


We have $x_1x_j, x_jx_7 \in E(G_1)$. Since $3 \leq j \leq 4$, $x_1x_j$ and $x_jx_7$ are diagonals of $G$. By Lemma~\ref{diagonallemma}, $x_1x_j$ partitions $G$ into two mops $G_{1,1}$ and $G_{1,2}$ with $V(G_{1,1}) = \{x_1, x_2, \dots, x_j\}$, and $x_jx_7$ partitions $G$ into two mops $G_{2,1}$ and $G_{2,2}$ with $V(G_{2,1}) = \{x_j, x_{j+1}, \dots, x_7\}$. By Lemma~\ref{maxofn2}, at most one of $x_1$ and $x_j$ is of degree $2$ in $G_{1,1}$, and at most one of $x_j$ and $x_7$ is of degree $2$ in $G_{2,1}$. Thus, by Lemma~\ref{basicfacts}(b), $d_{G_{1,1}}(y) = 2$ for some $y \in V(G_{1,1}) \backslash \{x_1,x_j\}$, and $d_{G_{2,1}}(z) = 2$ for some $z \in V(G_{2,1}) \backslash \{x_j,x_7\}$. We have $d_G(y) = d_{G_{1,1}}(y) = 2$ and $d_G(z) = d_{G_{2,1}}(z) = 2$. Also, $y, z \notin V(G')$. Now $d_{G}(v) = d_{G'}(v)$ for each $v \in V(G') \backslash \{x_1, x_7\}$, and, by Lemma~\ref{basicfacts}(b), at most one of $x_1$ and $x_7$ is of degree~$2$ in $G'$. Therefore, $n_2' \le n_2-1$.


We have $n'+n_2' \leq n+n_2-6$. By the induction hypothesis, $\ii_1(G') \le (n' + n_2')/6 = (n+n_2)/6 - 1$. Let $S'$ be a smallest $K_{1,2}$-isolating set of $G'$. Then, $|S'| = \ii_1(G') \le (n+n_2)/6 - 1$.  Since $x_j$ is adjacent to $x_1$ and $x_7$, $S' \cup \{x_j\}$ is a $K_{1,2}$-isolating set of $G$, and hence $\ii_1(G) \le |S'| + 1 \le (n+n_2)/6$.~\hfill{\smallqed}

\begin{claim}
\label{c2:3}
If $\ell = 7$, then $\ii_1(G) \le \frac{n+n_2}{6}$.
\end{claim}
\proof As in the proof of Claim~\ref{c1:3}, we may assume that $j = 4$.  Let $G' = G_2$ and $n' = |V(G_2)|$. We have $V(G') = V(G) \setminus \{x_2,x_3,x_4,x_5,x_6,x_7\}$, so $n' = n-6 \ge 5$. By the same reasoning used for Claim~\ref{c2:2}, there are at least two vertices in $V(G) \backslash V(G')$ of degree $2$ in $G$, at most one of $x_1$ and $x_8$ is of degree~$2$ in $G'$, and $d_{G}(v) = d_{G'}(v)$ for each $v \in V(G') \backslash \{x_1, x_8\}$. Thus, $n_2' \le n_2-1$, and hence $n'+n_2' \le n+n_2-7$. By the induction hypothesis, $\ii_1(G') \le (n'+n_2')/6 < (n+n_2)/6-1$. Let $S'$ be a smallest $K_{1,2}$-isolating set of $G'$. Then, since $S' \cup \{x_4\}$ is a $K_{1,2}$-isolating set of $G$, we have $\ii_1(G) \le |S'| + 1 = \ii_1(G') + 1 < (n+n_2)/6$.~\hfill{\smallqed}

By Claims~\ref{c2:1}, \ref{c2:2} and \ref{c2:3}, we may assume that $\ell = 8$. As in the proof of Theorem~\ref{t:main1} for the same case, we may assume that $j = 5$. Let $G' = G_2$ and $n' = |V(G_2)|$. We have $V(G') = V(G) \setminus \{x_2,x_3,x_4,x_5,x_6,x_7,x_8\}$, so $n' = n-7 \ge 4$. By the same reasoning used for Claim~\ref{c2:2}, there are at least two vertices in $V(G) \backslash V(G')$ of degree $2$ in $G$, at most one of $x_1$ and $x_9$ is of degree~$2$ in $G'$, and $d_{G}(v) = d_{G'}(v)$ for each $v \in V(G') \backslash \{x_1, x_9\}$. Thus, $n_2' \le n_2-1$, and hence $n'+n_2' \le n+n_2-8$. By the induction hypothesis (recall that the first bound in the theorem is being proved for $n \geq 4$), $\ii_1(G') \le (n'+n_2')/6 < (n+n_2)/6-1$. Let $S'$ be a smallest $K_{1,2}$-isolating set of $G'$. Then, since $S' \cup \{x_5\}$ is a $K_{1,2}$-isolating set of $G$, we have $\ii_1(G) \le |S'| + 1 = \ii_1(G') + 1 < (n+n_2)/6$.~\hfill{\qed}


\end{document}